\documentclass[letterpaper, 10pt, journal, final]{IEEEtran}
 \IEEEoverridecommandlockouts
\usepackage{cite} 
\usepackage{amsmath,amssymb,amsfonts} 
\usepackage{algorithmic} 
\usepackage[pdftex]{graphicx} 
\usepackage{textcomp} 
\usepackage{xcolor} 
\usepackage{multicol}
\usepackage{booktabs}
\usepackage{subfigure}
\usepackage{url}
\usepackage{breakurl}
\usepackage{soul} 
\usepackage[ruled,commentsnumbered]{algorithm2e}

\definecolor{gre}{rgb}{0.01, 0.75, 0.24}

\newcommand{\Hquad}{\hspace{0.3em}} 
 
\begin{document}
\title{The Added Value of Coordinating Inverter Control}
\author{Peter~Lusis,
        Lachlan~L.~H.~Andrew,\\
        Ariel~Liebman,
        and~Guido~Tack}
\maketitle 
\thispagestyle{plain} 
\pagestyle{plain} 
\begin{abstract} 
Coordinated photovoltaic inverter control with centralized coordination of curtailment can increase the amount of energy sent from low-voltage (LV) distribution networks to the grid while respecting voltage constraints. First, this paper quantifies the improvement of such an approach relative to autonomous droop control, in terms of PV curtailment and line losses in balanced networks. It then extends the coordinated inverter control to unbalanced distribution networks. Finally, it formulates a control algorithm for different objectives such as the fairer distribution of PV curtailment and rewarding PV customers for utilizing the excess power locally. 

The coordinated inverter control algorithm is tested on the 114-node and 906-bus LV European test feeders with cable sizes between 50\;mm$^2$ and 240\;mm$^2$ and validated with reference to OpenDSS. The results demonstrate that coordinated inverter control is superior when applied to high impedance LV networks and LV networks constrained by the distribution transformer capacity limits compared to autonomous inverters. On the 95\;mm$^2$ overhead line, it yields a 2\% increase on average in the utilized PV output with up to 5\% increase for some PV locations at higher penetration levels. Up to a 20\% increase in PV hosting capacity was observed for location scenarios with PV system clustering. 
\end{abstract}

\section*{Nomenclature}
\addcontentsline{toc}{section}{Nomenclature}
\begin{IEEEdescription}
[\IEEEsetlabelwidth{$\textbf{V}^{\text{min}},\textbf{V}^{\text{DB}}$}]
\item[$\mathcal{A}$] set of nodes $a$ with autonomous inverters.
\item[$\mathcal{C}$] set of nodes $c$ with coordinated inverters.
\item[$\mathcal{J}_k$] set of nodes $j$ with solar PV systems.
\item[$\mathcal{L}$] set of nodes $l$ with legacy inverters.
\item[$\mathcal{N}$] set of nodes $n$ with customers.
\item[$\mathcal{M}$] set of lines.
\item[$\mathcal{S}$] set of PV location scenarios $k$.
\item[$\textbf{P}^{\text{av}}$] available PV output.
\item[$\textbf{P}^{\text{d}}_n$] customer $n$ active power demand.
\item[$\textbf{Q}^{\text{d}}_n$] customer $n$ reactive power demand.
\item[$\textbf{Q}^{\text{min},pu}$] inverter reactive power support limit.
\item[$\textbf{R}$] line resistance.
\item[$\textbf{S}$] inverter kVA capacity.
\item[$\textbf{X}$] line reactance.
\item[$\textbf{y}^*$] complex conjugate of line admittance.
\item[$P_j$]  active power curtailment of inverter $j$.
\item[$P^{\text{inj}}_j$]  injected active power of inverter $j$.
\item[$P^{\text{util}}_j$]  active power curtailment of inverter $j$.
\item[$Q_j$] reactive power support of inverter $j$.
\item[$V_n$] voltage at node $n$ derived through optimization.
\item[$\hat V_n$] measured node $n$ voltage (wrt. OpenDSS).
\item[$V^{\text{max}}_j$]  maximum permitted voltage of inverter $j$.
\item[$V^{\text{nom}}_n$] nominal voltage at node $n$ about which linearization is performed.
\item[$\textbf{V}^{\text{trip}}$] 10-min average trip voltage.
\item[\smash{
\begin{IEEEeqnarraybox*}[][t]{l}
    \textbf{V}^{\text{min}}, \textbf{V}^{\text{DB}} \\
    \hphantom{V_1,{}} \textbf{V}^{\text{Qmin}}
\end{IEEEeqnarraybox*}
}] inverter droop curve voltage setpoints.\\ \mbox{}
\item[$\eta$] damping factor.
\item[$\theta$] the angle between the old and new reference frames.
\item[$\omega^\uparrow$] weights for autonomous inverter disconnection.
\item[$\omega^\downarrow$] weights for autonomous inverter reconnection.
\end{IEEEdescription}
\section{Introduction}\label{Introduction}
The first generation of inverters of photo-voltaic (PV) systems were designed to maximize active power injection into the grid. Such PV inverters, which we will call \emph{legacy} inverters, led to voltage quality problems in many low-voltage (LV) distribution networks. The risk of greater adverse effects during high PV generation periods has been addressed by adopting new inverter standards with grid support functions (GSF). However, more power quality and network-operational challenges arise with the continual increase in solar PV deployment. This trend has resulted in low-voltage areas where distribution network service providers (DNSPs) have to reject new solar PV connections. 

Voltage regulation methods, such as off-load tap changers, capacitors, and voltage regulators, were designed for one-way power flow \cite{Agalgaonkar2014, Singhal2018, Jothibasu2016}, and lack the fast switching ability required to follow the changes in PV generation output \cite{Watson2016}. Upgrading transformers and power lines is capital-intensive, and replacing power transformers will not solve all of the voltage issues on LV residential networks or long radial feeders \cite{Farivar2011}. An emerging alternative is to use PV inverters as a non-network alternative for mitigating operational challenges. That is the approach taken here. 

\emph{Overvoltage} is when voltage levels exceed statutory limits posing a risk to damage electric appliances. The first generation of inverters with grid support functions (GSF) use \emph{autonomous} voltage regulation following a droop curve, which relies solely on local voltage measurements. In the preferred form (Volt/VAr), inverters absorb or inject reactive power \cite{Farivar2012, Anonas2018, Seuss2016}. Increasing inverter reactive power support close to the end of the radial network improves the voltage profile more than fixed reactive power compensation  that ignores the size and location of PV inverters \cite{Ali2018}. However, we will use homogeneous support because consumers may consider heterogeneous solutions unfair.

Inverters can also control active power (Volt/Watt control). A method to optimally design Volt/Watt droop curve reference points for each inverter was presented in \cite{GhapandarKashani2017}, \cite{Kraiczy2018}. Inverter active power setpoints for equal power curtailment using Volt/Watt droop curve was implemented in \cite{Ali2015}, while \cite{Seuss2017} formulated proportional curtailment with a smoothing term to account for the action of uncontrollable voltage regulators. Both algorithms require manual parameter tuning. 

Volt/VAr and Volt/Watt control were combined in \cite{Giraldez2017} showing a lower overall PV curtailment in all Hawaiian Electric secondary circuits than either by itself. In high PV penetration scenarios in \cite{Giraldez2018} with Volt/VAr-Volt/Watt yields 0.3-0.63\% more generation than with Volt/VAr alone. However, \cite{Giraldez2018} only considered the networks with a small number of customers per distribution transformer. Droop control was shown to be less effective in suburban areas with a large number of customers per distribution transformer in \cite{Parajeles2017}. 

Decentralized voltage regulation typically leads to suboptimal operation of the network \cite{Baker2018, DallAnese2018}. This raises the question of whether the coordination of the power electronic equipment can provide significant benefit in managing distribution networks in the future. On the one hand, adding coordination increases the costs and complexity of the system. On the other hand, the communication equipment can be shared with tasks such as demand response and coordinated electric vehicle charging, and so need not be an additional expense \cite{Chapman2018}. 

\emph{Coordinated} inverter control involves collecting local voltage, load, and PV output measurements at a central location to calculate inverter setpoints such that network constraints and the operator's objective are met, and communicates the setpoints back to PV inverters. Improved voltage control through optimally controlling the power dispatch from inverters in response to renewable generation outputs is presented in \cite{DallAnese2014}. With the optimal inverter dispatch (OID) formulation, active and reactive power setpoints ($P$,$Q$) are updated continuously with the objective of minimizing power curtailment and line losses. Solar PV systems with microinverters can provide additional value as output from each panel can be optimized independently to maximize the system DC output \cite{Gagrica2015}. However, such an approach increases complexity significantly. 

OID has been formulated to work with clouds \cite{Ding2017},  uncertainty in load and PV \cite{Baker2018}, in the combination with other inverter control methods \cite{Andrew2020}, dynamics in the voltage \cite{Liu2018dynamic} and a combination of infrequent central control and fast local regulation \cite{Wang2019}. A combination of local voltage regulation with a periodical update of droop curve reference points for unbalanced three-phase four-wire networks was demonstrated in \cite{Weckx2014}, while a hierarchical droop-control model for voltage regulation with reduced communication between the supervisory node and local inverters was demonstrated in \cite{Karthikeyan2017}. However, the use of any kind of droop curve restricts the operational flexibility to a linear function compared to a case where inverters use the entire feasible operating region as shown in \cite{DallAnese2014, Ku2016, Lusis2019}.

Two inverter control strategies for mitigating sudden changes in load or solar PV output relying on a single-loop $P$-$Q$ regulation method and double loop control method were demonstrated in \cite{Jiandong2017}. Voltage regulation based on network clustering and two-stage inverter output optimization was presented in \cite{Chou2017}, but line losses were ignored.

Autonomous inverter control is the standard inverter-based voltage regulation method today, but, with increasing PV penetration, overvoltage disconnection will occur more often lowering the performance of PV systems, leaving negative voltage impacts on the upstream network stability.

In this paper, we formulate coordinated inverter control (CIC) operation designed to address voltage disturbances in PV-rich distribution networks during high solar irradiation periods, without sacrificing customers' ability to utilize PV output for their own demand. We then show that the biggest benefit of coordinated control over autonomous control is reducing line losses, and suggest that \emph{utilized} PV output, rather than curtailment, is the most appropriate metric for comparing PV inverter controls, particularly for networks with the smaller 95\;mm${}^2$ conductors. Finally, we show that coordinated inverter control can be applied to both balanced and unbalanced three phase networks, increasing PV hosting capacity and the utilized power output compared to autonomous and legacy inverters. The CIC algorithm is validated with reference to open-source software OpenDSS, which supports three-phase unbalanced network models \cite{OpenDSS2019}. The performance of coordinated and autonomous inverter control methods are compared by examining a 114-node LV circuit and IEEE 906-bus distribution network considering various line geometry to represent overhead and underground LV networks. 

The remainder of this paper is structured as follows. The inverter control strategies are defined in Section \ref{Formulation}, and the case study is described in Section \ref{CaseStudy}. The numerical results are presented and discussed in Section \ref{Results}. Final conclusions are provided in Section \ref{Conclusions}. 

\section{Inverter Control Algorithm}\label{Formulation}
The following notation will be used: all nodes are collected in the set $\mathcal{N}$, while $\mathcal{N}'\subset \mathcal{N}$ denotes the network without the slack bus (the secondary side of the distribution transformer). Complex voltages $V$ are expressed with $\Re$ and $\Im$ denoting the real and imaginary parts. Lines are represented as $(m,n) \in \mathcal{M} \subset \mathcal{N} \times \mathcal{N'}$. Subscripts $l\in\mathcal{L}, a\in\mathcal{A}$ and $c\in\mathcal{C}$ denote the quantities pertaining to each legacy inverter $l$, autonomous inverter $a$, and coordinated inverter $c$ at time $t$. System parameters beyond our control are written in bold. Other notation is as introduced throughout. 
\subsection{Coordinated Inverter Control}\label{OID}
Coordinated inverter control (CIC) gathers customers' inverter PV output data and voltages from the Advanced Metering Infrastructure (AMI) across the network to compute the optimum active power curtailment $P_c$ and reactive power support $Q_c$ for each inverter. These are then communicated back and implemented by the inverters\footnote{Current AMI updates readings comparatively slowly due to communication network limitations. In such circumstances, control should occur on two timescales. On the fast timescale, droop control prevents over-voltage, and on the slower timescale, the parameters of the droop controller are set dynamically such that the operating point is controlled by coordinated control. The implications of this are a suitable topic for future research.}. The inverter optimization problem attempts to minimize overall system losses expressed as power curtailment from coordinated inverters $\phi(P_{\mathcal{C},t})$ and line losses $\rho(V_t)$
\begin{align}
\underset{V,P_{\text{c}}}{\text{minimize \textit{}}} 
& \phi(P_{\mathcal{C},t}) + \rho(V_t) 
\label{eq:Objmain} \\
\text{subject to } & \eqref{eq:InvConst1}-\eqref{eq:VmaxOID}. \nonumber
\end{align} 
The active power curtailment $P_{c}$ for the set of coordinated inverters $\mathcal{C}$ is
\begin{equation} 
\phi(P_{\mathcal{C},t})= \sum_{c\in\mathcal{C}} P_{c,t}.
\label{eq:Obj2} 
\end{equation} 
Line active power losses are 
\begin{align}
\rho(V_t) = \sum_{(m,n)\in\mathcal{M}}&\Re\{\textbf{y}^*_{nm}\}\Big((\Re\{V_{m,t}\}-\Re\{V_{n,t}\})^2 \nonumber \\
& +(\Im\{V_{m,t}\}-\Im\{V_{n,t}\})^2\Big), 
\label{eq:Obj1} 
\end{align}
where $\textbf{y}^*_{mn}$ denotes the complex conjugate of admittance between nodes $m$ and $n$. PV curtailment $P_c$ is limited by the excess power, defined as the difference between the available active power on the AC side of the inverter $\textbf{P}^\text{av}$ and customer load $\textbf{P}^\text{d}$:
\begin{align}
  P_{c,t} = \begin{cases}
    0 \leq P_{c,t} \leq \textbf{P}^\text{av}_{c,t}-\textbf{P}^\text{d}_{c,t}, & \text{if } \textbf{P}^\text{av}_{c,t}>\textbf{P}^\text{d}_{c,t},\quad \forall c \in \mathcal{C},\\
     0, & \text{otherwise}.
  \end{cases}
 \label{eq:InvConst1}
\end{align}
The active power injection and reactive power support are bounded by the inverter's rated apparent power limit  
 \begin{equation}
(\textbf{P}^\text{av}_{c,t}-P_{c,t})^2 + (Q_{c,t})^2 \leq \textbf{S}^2, \quad \forall c \in \mathcal{C}.
\label{eq:InvConst2}
\end{equation}

Since our focus is entirely on the inverter operation during high voltage periods, we ignore evening and cloudy periods and operate inverters only in the lagging mode. Then, reactive power support is limited to $
\textbf{Q}^{\text{min},pu}\textbf{S} \leq Q_c \leq \textbf{0}
$. 

The nodal voltage balance is based on a linear approximation of power flow equations introduced in~\cite{Dhople2016}. This approach has been previously applied in~\cite{Guggilam2016}, demonstrating fast convergence properties.  In our case, it is\footnote{Better results can be obtained by dividing these by $V^{\text{nom}}$ to give a first order Taylor approximation.}
\begin{align}
&\Re\{\Delta V_{n,t}\} = 
\sum_{c\in\mathcal{C}}\Big( \textbf{R}_{nc}(\textbf{P}^\text{av}_{c,t}-P_{c,t}) + \textbf{X}_{nc}Q_{c,t}\Big) \nonumber \\ &
-\sum_{m\in\mathcal{\mathcal{N}}}\Big(\textbf{X}_{nm}\textbf{Q}^\text{d}_{m,t}  
+ \textbf{R}_{nm}\textbf{P}^\text{d}_{m,t}\Big) 
\quad \forall n\in \mathcal{N},
\label{eq:NodalConst1}
\end{align}
\begin{align}
& \Im\{\Delta V_{n,t}\} = \sum_{c\in\mathcal{C}}\Big( \textbf{X}_{nc}(\textbf{P}^\text{av}_{c,t}-P_{c,t})  -  \textbf{R}_{nc}Q_{c,t}\Big)  \nonumber \\ &
+\sum_{m\in\mathcal{N}} \Big(\textbf{R}_{nm}\textbf{Q}^\text{d}_{m,t} - \textbf{X}_{nm}\textbf{P}^\text{d}_{m,t}\Big), \quad \forall n\in \mathcal{N}.
\label{eq:NodalConst2} 
\end{align}

$\mathbf{R}_{nm}$ and $\mathbf{X}_{nm}$ denote the real and imaginary parts of the line impedance matrix (the inverse of the bus admittance matrix).  We assume that line impedances and network topology are known when solving the quadratically constrained quadratic program (QCQP). Active and reactive loads are denoted by $\mathbf{P}^\text{d}$ and $\mathbf{Q}^\text{d}$. 

The maximum voltage is limited by $ \Re\{V_{c,t}\} = V^\text{nom} + \Re\{\Delta V_{c,t}\} \leq V_{c,t}^\text{max}$. The accuracy of the linearization is best when the perturbations are small.  For this reason, we do not linearize around 1 per unit, but instead adjust the nominal point to match the measured operating point,
\begin{equation}
    V^\text{nom}_{n,t} = V^\text{nom}_{n,t-1} + \eta(\hat{V}_{n,t-1} - V_{n,t-1}), \quad \forall n \in \mathcal{N},
    \label{eq:Vnom}
\end{equation}
where $\hat{V}_{n,t-1}$ is the voltage measurement at the previous time interval received from the AMI\footnote{In this study, obtained from OpenDSS.} and the model-calculated voltages on each node. The damping factor and $\eta$ = 0.4 is a step size is derived empirically to generate the closest solution to the observed voltage values.

Coordinated inverters are operated below the maximum voltage threshold $V_{c,t}^\text{max}$, which corresponds to the average tripping voltage $\mathbf{V}^\text{trip}$ of conventional droop-based (autonomous) inverters. Since CIC relies on bi-directional communication, fast-changing network conditions may sometimes cause the actual nodal voltages to exceed the calculated CIC voltages the inverter target values are set for. Moreover, network interruption or communication delays are also possible. This becomes particularly important when operating at voltages close to the tripping threshold $\textbf{V}^\text{trip}$. Therefore, we add a local constraint to update the maximum voltage threshold such that
\begin{align}
  V_{c,t}^\text{max} = \begin{cases}
    V_{c,t}^\text{max} - (\hat{V}_{c,t-1}- \textbf{V}^\text{trip})
    , & \text{if } \hat{V}_{c,t-1}>\textbf{V}^\text{trip} \\
     V_{c,t}^\text{max} , & \text{otherwise}.
  \end{cases}
 \label{eq:VmaxOID}
\end{align}
\subsection{Applying CIC to unbalanced distribution networks}
Load imbalance is often observed in three-phase four-wire distribution networks. We extend the single phase equivalent model to account for line coupling effect between phases, following the formulation demonstrated in \cite{Koirala2019}. The phase impedance matrix $\textbf{Z}$ is a symmetrical square matrix with line self-impedance $\textbf{R}^{ii}_{nm} + j\textbf{X}^{ii}_{nm}$ on the principal diagonal and mutual impedance $\textbf{R}^{ij}_{nm} + j\textbf{X}^{ij}_{nm}$ on the off-diagonal, where $i,j\in\{A,B,C\}$.

A three-phase network model can be described using voltage phasors that are separated by 120$^\circ$ and rotating in a three-phase domain. The injections of active power ($p = \textbf{P}^\text{av}-\textbf{P}^d-P$) and reactive power ($q = Q-\textbf{Q}^d$) initiate voltage changes on each node. We apply the Park transformation to project the real and imaginary voltage perturbations $\Delta V$ around $\textbf{V}^\text{nom}$ from the Cartesian reference frame onto a new reference frame
\begin{align} 
\begin{bmatrix} 
\Re\{\Delta V_{nm,i,t}\} \\
\Im\{\Delta V_{nm,i,t}\} \\
\end{bmatrix}
=
\textbf{D}
\begin{bmatrix} 
 \Big(\textbf{R}^{ij}_{nm}p_{m,i,t}+\textbf{X}^{ij}_{nm}q_{m,i,t}\Big)\\ 
 \Big(\textbf{X}^{ij}_{nm}p_{m,i,t}-\textbf{R}^{ij}_{nm}q_{m,i,t}\Big)\\
\end{bmatrix}, \\
\Hquad i\in\{A,B,C\}, \Hquad \forall m,n \in \mathcal{N},\nonumber
\end{align}
where $\textbf{D}$ is the Park transform matrix 
\begin{equation}
\textbf{D}=
\begin{bmatrix} 
\cos\theta & \sin\theta \\ 
-\sin\theta & \cos\theta\\ 
\end{bmatrix},
\end{equation}
and $\theta \in \{0, 120^\circ, -120^\circ\}$ is the angle between the rotating reference frame and the Cartesian reference frame. From here we can calculate nodal balance equations as 
\begin{align}
\Re\{V_{n,i,t}\} = \textbf{V}^\text{nom}_{n,i,t} +\sum_{m\in\mathcal{N}}\Re\{\Delta V_{nm,i,t}\},~i \in\{A,B,C\},~\forall n \in \mathcal{N},\nonumber\\
\Im\{V_{n,i,t}\} = 
\sum_{m\in\mathcal{N}}\Im\{\Delta V_{nm,i,t}\},~i \in\{A,B,C\},~\forall n \in \mathcal{N}.
\end{align}

Voltage magnitude is recovered as
\begin{equation}
    |V_{n,i,t}| = \sqrt{\Re\{V_{n,i,t}\}^2+\Im\{V_{n,i,t}\}^2},~ i\in\{A,B,C\},~\forall n \in \mathcal{N}.
\end{equation}
The angular orientation of voltage magnitude in reference to the real axis satisfies $\tan(\angle V) = I / R$. Numerical results of coordinated inverter control application for unbalanced distribution networks are validated with reference to OpenDSS (Table~\ref{tab:valid}). The pseudo-code for coordinated inverter algorithm is given as follows
\begin{algorithm} 
\caption{Coordinated Inverter Control}\label{algorithm_MixOID_full}
\SetAlFnt{\small\sf} 
\SetKwFunction{Union}{Union}\SetKwFunction{FindCompress}{FindCompress}
\SetKwInOut{Input}{Input}\SetKwInOut{Output}{Output}
\Input{Given initial $V^{\text{nom}}_{n,0}=\textbf{1}, \hat V_{n,0}=\textbf{1}$; 
}
initialisation\;
\For{each time $t$}{
Collect PV output $\textbf{P}^{\text{av}}_{c,t}$ and load data $\textbf{P}^{d}_{n,t}, ~  \textbf{Q}^{d}_{n,t}$\;
\For{$c\in\mathcal{C}$}{
Solve optimization problem \eqref{eq:Objmain} to find {$P_{c,t}, ~ Q_{c,t}, V_{n,t}$}
}
Measure voltages $\hat V_{n,t}$ and update point about which linearization $V^{\text{nom}}_{n,t+1}$ is performed according to \eqref{eq:Vnom}\;
\If{$\hat V_{c,t}>\textbf{V}^{\text{trip}}$}{
Update the maximum permitted voltage threshold $V^{\text{max}}_{c,t+1}$ according to \eqref{eq:VmaxOID}\;
}
}
\end{algorithm}
\subsection{Autonomous Inverter Control}\label{PassiveInverter}
Autonomous inverter control is network-agnostic, in that each autonomous inverter  $a\in \mathcal{A} \subseteq \mathcal{N'}$ utilizes only local voltage measurements to determine operational setpoints. We use a combination of Volt/VAr and Volt/Watt droop curves. The active zone of each droop curve is a linear function with fixed reference points. Since the goal is to maximize the use of the available active power  $\textbf{P}^{av}$, inverters are operated in reactive power priority mode, hence voltage regulation is first attempted through reactive power support. The nominal inverter reactive power injection in the grid operating in Volt/VAr response mode is 
\begin{equation}
Q^{pu}_{a,t} = \textbf{m}^\text{V}V_{a,t}+\textbf{c}^\text{V},\quad \forall a \in \mathcal{A}.
\label{eq:Q_VV}
\end{equation}
The slope $\textbf{m}^\text{V}$ and intercept $\textbf{c}^\text{V}$ of the Volt/VAr droop curve are
\begin{align} 
\textbf{m}^\text{V}=\textbf{Q}^{\text{min},pu}/({\textbf{V}^{\text{Qmin}}-\textbf{V}^\text{DB}}),
\label{eq:slopeVV2} 
\\
\textbf{c}^\text{V}=\textbf{Q}^{\text{min},pu}\textbf{V}^\text{DB}/(\textbf{V}^{\text{Qmin}}-\textbf{V}^\text{DB}), 
\label{eq:intercVV2} 
\end{align}
where $\textbf{V}^\text{DB}$ is the upper reference point of the deadband, and $\textbf{V}^{\text{Qmin}}$ denotes the voltage level at which the maximum allowable reactive power support is reached. When operating in Volt/Watt response mode, $Q^{pu}_{a,t} = \textbf{Q}^{\text{min},pu}$. Converting to base units, this gives $Q_{a,t} = Q^{pu}_{a,t}\textbf{S}_au_{a,t}$. The inverter status is determined by $u_{a,t}\in\{0,1\}$. In the simulations all inverters are initially set to ON ($u_{a,0}=1$).

When voltage levels cannot be maintained within operational limits using reactive power, the Volt/Watt function is activated, linearly reducing active power output until the cut-off voltage $\textbf{V}_a^{\text{max}}$ is reached:
\begin{equation}
P^{pu}_{a,t}= \textbf{m}^\text{W}V_{a,t}+\textbf{c}^\text{W}.
\label{eq:Q_VW}
\end{equation}
The slope $\textbf{m}^\text{W}$ and intercept $\textbf{c}^\text{W}$ are calculated as a linear function passing through points ($\textbf{V}^\text{Qmin}, \textbf{P}^\text{av}$) and ($\textbf{V}^{\text{max}}, \textbf{P}^{\text{min}}$). Then, the inverter injected power $P^\text{inj}_{a}$ operating in the Volt/Watt mode can be calculated as 
\begin{equation}
  P^\text{inj}_{a,t} = \min\left\{\sqrt{\textbf{S}_a^2-Q_{a,t}^2},\textbf{P}^\text{av}_{a,t}P^{pu}_{a,t}\right\}u_{a,t}.
 \label{eq:Pinj_VW}
\end{equation}

\begin{figure} 
\centerline{\includegraphics[width= 9 cm]{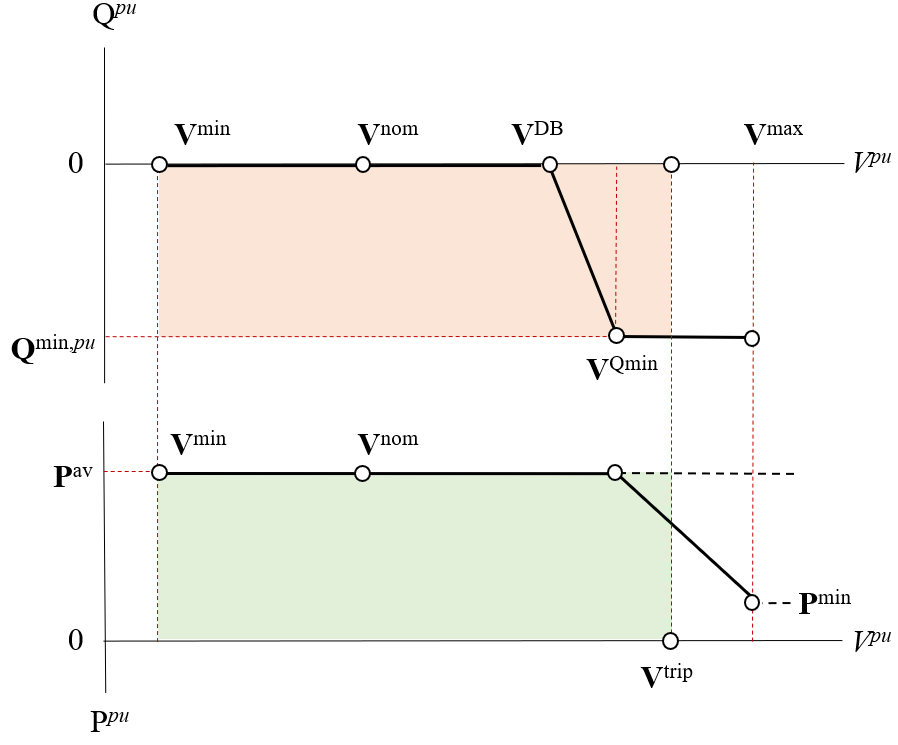}}
\caption{Volt/VAr and Volt/Watt droop curves. The highlighted area represents the feasible setpoint space for coordinated inverter control.} 
\label{fig:droopcruves} 
\end{figure}

Autonomous inverters have multiple conditions that govern their disconnection from the grid. If none is satisfied, the inverters will follow Volt/VAr or Volt/Watt response modes.  If the moving average voltage for $\textbf{u}^{\uparrow}$ time intervals is above the \textit{average trip} voltage $\textbf{V}^\text{trip}$, autonomous inverters will shut off. 

In our simulations, which use discrete time steps of 1~minute, the change in demand and available PV output between consecutive periods can lead to multiple inverters exceeding $\textbf{V}^\text{trip}$ and shutting off. The concurrent adjustment of inverter setpoints may lead to unwanted voltage oscillations as demonstrated in \cite{Jahangiri2013}. The fact that this is precisely simultaneous is merely an artefact of the simulation. To avoid this artefact, we allow only one autonomous inverter $a$ to disconnect in each period based on random weighted sampling with larger weights $w^{\uparrow}_{a,t}$ assigned to nodes with higher voltages
\begin{equation}
w^{\uparrow}_{a,t} =  (\textbf{V}^{\text{max}}-V_{a,t})^{-2}.
\label{eq:weightedSamplingDisc}
\end{equation} 
Without this constraint, the simulations amplify voltage oscillations beyond what would be found in a real network. Capacitor banks with oscillation damping capabilities or other power electronics can address this issue in the actual network while the one-at-a-time rule is needed only for simulations.
\par
The average trip voltage rule is overridden if $V_{a,t}\geq\textbf{V}_a^{\text{max}}$, which leads to instantaneous disconnection of all inverters $a$ that meet this condition. Inverters remain disconnected for at least $\textbf{u}^{\downarrow}$ periods, and will reconnect if the voltage has returned below $\textbf{V}^{\text{trip}}$. Again, only one simulated inverter can reconnect at a time based on weighted random sampling with weights 
\begin{equation}
w^{\downarrow}_{a,t} = (V_{a,t}-\textbf{V}^\text{nom})^{-2}.
\label{eq:weightedSamplingRec}
\end{equation} 
\subsection{Legacy Inverters}\label{LegacyInverter}
Legacy inverters located on nodes $l\in\mathcal{L} \subset \mathcal{N'}$ can operate only at unity power factor. During the periods with high voltage observations, the tripping conditions of legacy inverters and the inverter ON/OFF status $u_{l,t}$ are governed by  $\textbf{V}^\text{\text{trip}}$ and $\textbf{V}_l^{\max}$ similar to autonomous inverters (see Appendix). The injected active power $P_l^\text{inj}$ from legacy inverters is 
 \begin{equation}
   P^\text{inj}_{l,t}=\begin{cases}
     \textbf{P}^\text{av}_{l,t}, & \text{if } u_{l,t}=1,\quad \forall l \in \mathcal{L},\\
     0, & \text{otherwise}.
   \end{cases}
 \label{eq:LegacyInv}
 \end{equation}
Time-series power flow analysis for autonomous and passive inverter control is implemented in OpenDSS software and operated via the COM interface through Matlab.
\section{Analysis Framework}\label{CaseStudy}
\subsection{Evaluation Criteria \& Scenario Generation} 
This paper aims to quantify the added benefits from coordinated inverter control in terms of the number of systems that can be accommodated, energy fed back to the grid, transient overvoltage, and fairness between payments to customers at different locations. The likelihood of experiencing operational issues due to the presence of solar PV depends on the aggregated capacity of solar PV and the location of each  PV system \cite{Kim2018}, \cite{Rylander2016NREL}. We characterize the number of systems that can be supported by a pair of values. The lower bound of PV hosting capacity $cap^{\min}$ is defined as the lowest PV customer penetration level such that in at least in one scenario some customers experiences PV curtailment. The upper bound $cap^{\max}$ is the smallest penetration level such that curtailment is inevitable in all scenarios. That is, 
\begin{align} 
cap^{\min}= {\text{min \textit{}}} \{x: (\exists~k\in\mathcal{S}: |\mathcal{J}_k|/|N| \le x),\exists~t: P_{j,t} > 0\}
\label{eq:hosting_min}  
\\
cap^{\max}=
{\text{min \textit{}}} \{x: (\forall~k\in\mathcal{S}: |\mathcal{J}_k|/|N| \ge x),\exists~t: P_{j,t} > 0\} 
\label{eq:hosting_max}  
\end{align}
where the location scenarios constitute the set $\mathcal{S}$, and in each case $k\in \mathcal{S}$, $\mathcal{J}_k\subseteq\mathcal{N}$ is the set of nodes $j$ with solar PV systems (for all inverter control methods). 

To capture the effect of non-uniformity in the PV penetration, two PV deployment scenarios represent the cases where PV systems are clustered around the points of the minimum and maximum effective impedance $Z$ from the distribution transformer. Another 18 PV deployment scenarios were generated at each customer penetration level by applying repeated random sampling to simulate customers independently obtaining PV systems. 

The utilized active power is 
\begin{align} 
\sum_{t\in\mathcal{T}}{P^{\text{util}}_{t}} = \sum_{t\in\mathcal{T}}\textbf{P}^{\text{av}}_t -
\phi(P_{j,t}) + \rho(V_t),
\label{eq:P_useful}  
\end{align} 
subtracting active power curtailment $\sum_{t\in\mathcal{T}}\phi(P_{j,t})$ and line losses $\sum_{t\in\mathcal{T}}\rho_t(V)$ from the available PV output. An argument for deploying coordinated inverter control would be a significant reduction in PV curtailment, line losses or both, to justify additional costs related to enabling and maintaining communication infrastructure between the inverters and the central node. The network operation for a summer day is simulated using each inverter control mode at 20 discrete customer penetration levels between 10 and 100\% in increments of 10\%. 
\subsection{Network Topology} 
The analysis described in the previous section is applied to the 114-node semi-urban low-voltage network of \cite{Prettico2016}, and IEEE 906-bus European test feeder with 55 loads \cite{Espinosa2015}, suggested for LV network studies by the IEEE PES Distribution Systems Analysis Subcommittee \cite{Schneider2017}. The size of solar PV systems is 6\,kWp, corresponding to an average rooftop PV system in Australia installed to date. The use of uniform nameplate capacity for all PV systems can be justified by implementing a number of random distributions of PV systems across the network. This was de-rated by about 17\% to account for reduced PV output below the rated capacity due to factors such as soiling of the panels, shading and aging, giving a maximum PV output of 5\,kW at the AC side of the inverter. The inverter reactive power limit is set to 0.44 lagging of the inverter kW output. The step size in~\eqref{eq:Vnom} was $\eta = 0.4$. 

Half-hourly household active power consumption from \cite{Ausgrid2012} was interpolated to 1-minute intervals using cubic splines. Data from the $i$th household, $i=1,\dots,30$ was allocated to buses with numbers $i\mod 30$. Simulations run from 8\,am to 7.30\,pm. Outside these hours, none of the inverter overvoltage response modes was triggered, even at 100\% penetration. The average load per household $\tilde{\textbf{P}}^\text{d}$  within this period is 0.77\,kW with variance of 0.27 \,kW. Reactive power demand $\textbf{Q}^\text{d}$ is set to $0.328~\textbf{P}^\text{d}$, giving a constant power factor of 0.95 leading. 

Autonomous inverter control settings are standardized within a DNSP network area, often across multiple states or jurisdictions. However, low-voltage networks vary  significantly with respect to the type of customers, topology, and other variables. Since coordinated inverter control entails additional implementation and operational costs, one would expect to deploy it when network hosting capacity or utilized power can be increased. 

\begin{table}[t]
\caption{Distribution line characteristics.}
\label{tab:cableModels} 
\centering
\resizebox{\columnwidth}{!}{%
\begin{tabular}{l c c c c c}
\toprule
LV line type & Area & 
\textbf{R}$^{ii}$ & \textbf{X}$^{ii}$ &
\textbf{R}$^{ij}$ & \textbf{X}$^{ij}$ \\
 & [$mm^2$]  & 
[$\Omega/km$] &[$\Omega/km$] &
[$\Omega/km$] & [$\Omega/km$] \\
\midrule 
Open wire \cite{Prettico2016} & 4x50 & 0.699 & 0.149 & 0.049 & 0.164 \\
Underground \cite{Koirala2019} & 4x70 & 0.759 & 0.243 & 0.316 & 0.193 \\
Open wire \cite{Safitri2015}  & 4x95 & 0.452 & 0.270 & 0.049 & 0.164 \\
Underground \cite{Weckx2014}& 4x150 & 0.227 & 0.078 & 0.070 & 0.078 \\
Underground \cite{Procopiou2019} & 3x240 & 0.072 & 0.199 & 0.021 & 0.048 \\
\bottomrule
\end{tabular}
}
\end{table}
The type and size of conductors installed in the low-voltage network may significantly affect the operation of coordinated and autonomous inverter control. To study this effect, we test autonomous and coordinated inverter control on five commonly installed low-voltage over-head and underground conductors. The cable specs with references are given in Table~\ref{tab:cableModels}.

\section{Simulation Results}\label{Results}
\subsection{Utilized PV output and hosting capacity}
\begin{figure}[tb] 
\centerline{\includegraphics[width=  \linewidth]{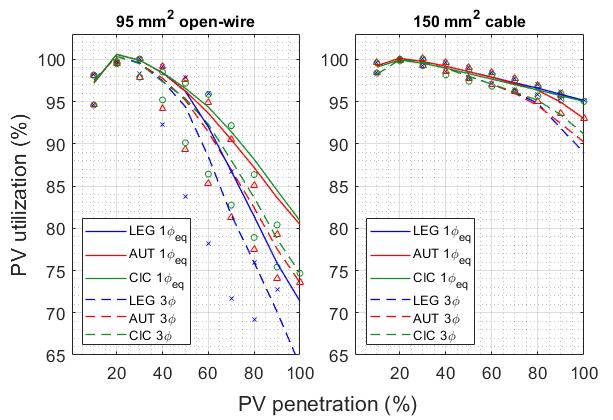}}
\vspace{0.1cm}
\centerline{\includegraphics[width= \linewidth]{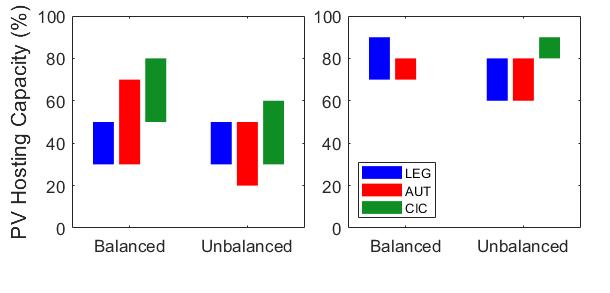}}
\caption{The total utilized PV output as a percentage of solar PV output proportional to PV penetration levels with 95mm$^2$ and 150m$^2$ conductors. The utilization rate accounts for PV curtailment and line losses that would occur compared to a scenario without solar PV systems.
Solid lines represent the average utilized PV output for a balanced network while dashed lines correspond to an unbalanced network. The top and bottom marker types denote two cases when PV are clustered at the end of the line or closest to the distribution transformer. The bottom figure shows PV hosting capacity for legacy (LEG), autonomous (AUT) and coordinated (CIC) inverters for the respective cables.}
\label{fig:fig2_PVutil}
\end{figure}

The first study varies customer PV penetration levels with uniformly sized conductors of 95\,mm${}^2$ and 150\,mm${}^2$ for each inverter control model applied on the 114-node LV network to compare utilized PV output levels (Fig.~\ref{fig:fig2_PVutil}). Smaller conductors of 50\,mm${}^2$ and 75\,mm${}^2$ were also considered. However, due to high impedance, those conductors breached the operating voltage range of $\{-6\%,+10\%\}$ even at low penetrations. Meanwhile, no voltage violations were observed with 240\,mm${}^2$ conductors. The difference between the cases when PV systems are clustered at the end of the line opposed to adjacent to the transformer reached 18\% for legacy inverters and 10\% for autonomous and coordinated inverters. This highlights the sensitivity of the results to the PV location in the network; resulting in a challenge for distribution system operators to determine a single value as the hosting capacity limit (Fig.~\ref{fig:fig2_PVutil}, bottom). Moreover, the larger 150\,mm${}^2$ cables prevented any curtailment with coordinated inverter control even at maximum penetration. The results suggest that setting a single hosting capacity limit for all low-voltage networks is a very conservative approach and may unnecessarily limit new connections. Therefore, PV hosting capacity should be determined based on the recorded conductor characteristics.

As the penetration level reaches 20\%, the utilized PV output is seen to surpass 100\%, demonstrating a reduction in line losses in comparison to a case without solar PV. In addition, at 20\% penetration level, if either all PV systems are located close to the transformer or all are at the end of the line (see the markers in Fig.~\ref{fig:fig2_PVutil}), the output utilization rate is below the average. This indicates that PV clustering leads to higher energy losses as the excess PV output has to be distributed to the parts of the network. 

In terms of utilized PV output, coordinated inverters consistently perform better than legacy inverters designed to maximize PV output injection into the grid. This margin diminishes with larger cables as fewer overvoltage events are observed. When compared to autonomous inverters, applying inverter coordination in unbalanced networks increased the PV utilization rate by 2\% on average, reaching 5\% for specific PV location scenarios. The benefit of coordinated inverter control is topology-dependent, so we distinguish between two cases.

First, voltage rise at lower penetration levels can be associated with PV clustering in a certain part of the network (often combined with phase unbalance). Such a situation leads to at least one autonomous inverter entering the Volt/Watt curtailment mode at 253\;V. This will consequently reduce the hosting capacity, which we define as the minimum PV penetration when curtailment occurs at any of the PV locations. In comparison, CIC calculates the role of each inverter in the reactive power provision to delay the first curtailment, while maintaining voltage below the tripping level (257\;V). As a result, we observed up to 20\% increase in PV hosting capacity.

\begin{figure}[b]
\centerline{\includegraphics[width= \linewidth]{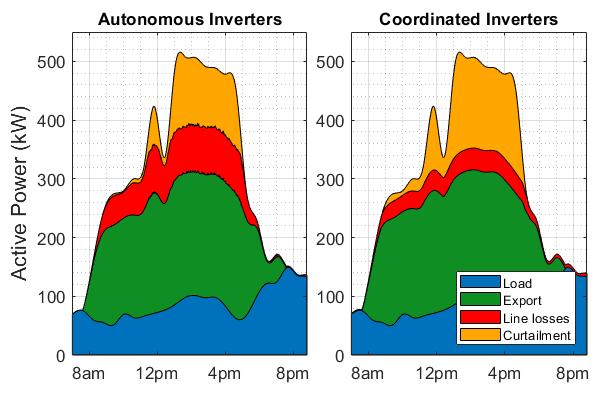}}
\caption{The breakdown of PV output by end-use with autonomous (on left) and coordinated (on right) inverter control. 90\% penetration level with 95\,mm${}^2$ conductors.} 
\label{fig:fig3_90penCompare}
\end{figure} 
Occasionally, autonomous inverters may yield lower overall PV curtailment than coordinated inverter control, however, overall line losses remain many times higher (Fig.~\ref{fig:fig3_90penCompare}). This is due to an increased reactive power demand when inverters operate in the Volt/VAr and Volt/Watt response mode. In contrast, coordinated inverter control explicitly maximizes PV output injected into the grid subject to the voltage constraints. The optimum setpoints returned to the coordinated inverters show that for some PV location scenarios lower overall energy losses can be achieved by using little to no reactive power support, as the additional energy injected by PV inverters is offset through increased line losses. Therefore, we recommend to use utilized PV power (local demand met plus export reaching the grid) instead of curtailment for the evaluation of inverter control methods.
\subsection{Distribution transformer loading}
\begin{figure}[tb]
\centerline{\includegraphics[width= \linewidth]{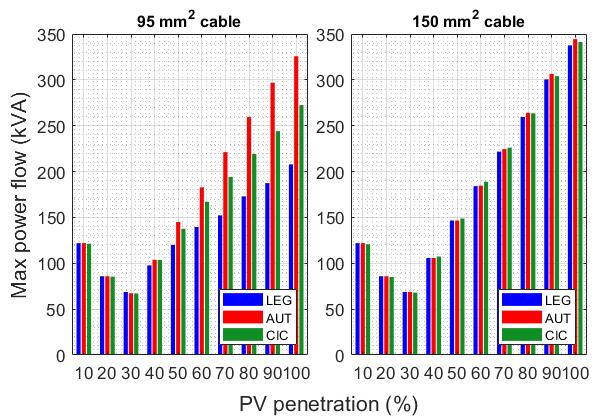}}
\caption{Maximum apparent power flow (kVA) through the distribution transformer recorded in any direction with 95\,mm$^2$ over-head and 150\,mm$^2$ underground conductors.} 
\label{fig:fig4_transformers}
\end{figure}

Another factor limiting PV hosting capacity is distribution transformer loading. In a network without solar PV systems, the highest apparent power occurs when the most power flows from the grid, and, as such, this peak decreases as PV penetration increases  (Fig.~\ref{fig:fig4_transformers}). However, from around 30\% penetration, the peak transformer loading occurs when PV output is exported to the grid. Note that the maximum apparent power does not dip to zero at the point when this reversal occurs, since the two maxima occur at different times in the day. 

As shown earlier, the implementation of coordinated inverter control in a low impedance network results in higher utilized power exports than autonomous inverters. Thus, the higher apparent power flow through the transformer with autonomous inverter control in Fig.~\ref{fig:fig4_transformers} is associated with an increase in reactive power demand. For example, if a low voltage network with 95\,mm$^2$ conductors was connected to a 300\,kVA transformer, the maximum hosting capacity to prevent overloading with autonomous inverters would be, at most, 90\%, while the 300\,kVA limit is not breached even at 100\% penetration of 6\,kWp PV systems under coordinated inverter control. Although reactive power support would lower voltages on the PV nodes, it is not used by CIC since the additional real power injection would be offset by higher line losses. The difference in transformer loading between coordinated and autonomous inverter scenarios reduces for 150\,mm$^2$ conductors, as lower impedances allow all power to be exported without excess reactive power demand.

\subsection{The impact of conductor size}
\begin{figure}[tb]
\centerline{\includegraphics[width= \linewidth]{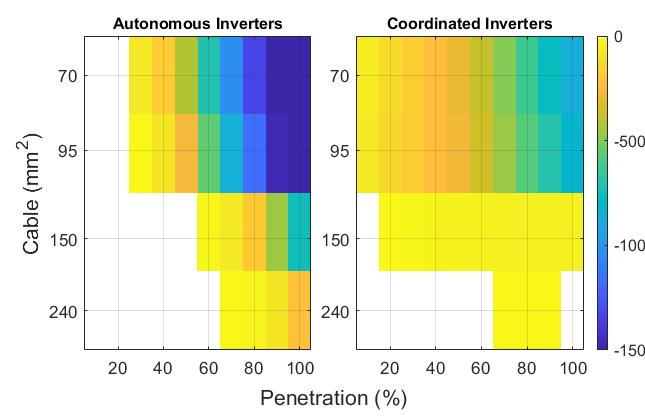}}
\caption{Average reactive power support (kVArh/day) over all PV locations for each conductor size and penetration level.} 
\label{fig:fig5_Qsupport}
\end{figure}

\begin{figure}[tb]
\centerline{\includegraphics[width= \linewidth]{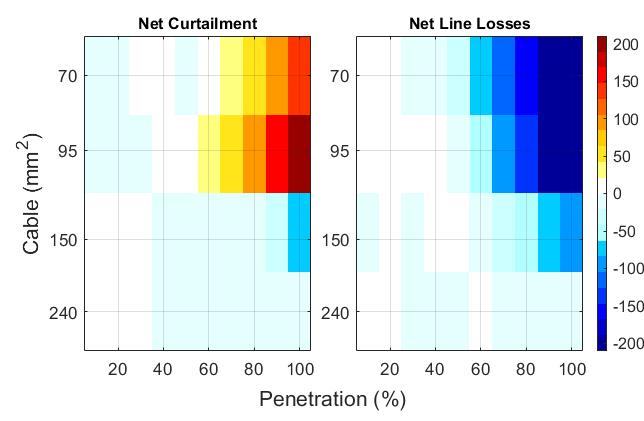}}
\caption{Average difference in curtailment (kWh/day) and line loss (kWh) between coordinated and autonomous inverter control over all PV locations for each conductor type and penetration level. The blue color (negative values) illustrates cases where coordinated inverter control perform better than autonomous inverter control. The opposite is true when the red color is used. Marginal difference in means is left as white space.}
\label{fig:fig6_net_Pcurt}
\end{figure}

The results of reactive power support with autonomous and coordinated inverter control for four different cable sizes are shown in Figs.~\ref{fig:fig5_Qsupport} and ~\ref{fig:fig6_net_Pcurt}. With autonomous inverter control, reactive power support increases with the cable size and penetration levels. Larger line impedance leads to overvoltages at much lower PV output levels, increasing the number of time periods over the day when inverters are required to operate in Volt/VAr and Volt/Watt response modes.

Coordinated inverted control demonstrates that at high penetration levels reactive power support is more effective when applied to the 95\,mm$^2$ conductors as opposed to higher gauge cables. As previously shown, smaller cables tend to contribute to higher line losses, and so reactive power support has limited capability to minimize the objective value. In contrast, large conductors do not cause overvoltages, avoiding the need for reactive power support. 

These results suggest that coordinated inverter control deployed on a network with small conductors could prevent the excess reactive power demand that occurs with autonomous inverters. When reactive power demand can be compensated using, for example, capacitors, or when the goal is to maximize the PV output injection into the grid, the line loss term can be removed from the objective function. This allows coordinated inverters to use ``freely'' available reactive power. Coordinated inverters also has the ability to utilize reactive power for voltage imbalance reduction across phases. This explains the use of reactive power for the larger cables in Fig.~\ref{fig:fig5_Qsupport} where reactive power is used even in absence of overvoltages as this leads to lower line losses.

Figure~\ref{fig:fig6_net_Pcurt} illustrates the difference in the amount of PV curtailment between autonomous and coordinated inverter control. It further supports the finding that the benefit is greatest for medium-sized cables as coordinated inverter control results in less PV curtailment for the customers. With the 95\,mm$^2$ conductors, coordinated inverter control curtails slightly less at penetration levels between 15-30\% due to a higher PV hosting capacity. At higher penetration levels, the difference in curtailment flips, and coordinated inverter control curtails more PV output than autonomous inverters in order to minimize line losses.

There was minimal difference in line losses with the 150 mm$^2$ cables, which shows the effect of higher reactive power demand with autonomous inverter control being offset by coordinated control exporting more active power. 
Taking a larger cable size leads to fewer voltage problems, thus the difference between control methods become less obvious, although coordinated inverter control still dominates at all penetration levels.
\subsection{Fairness}
\begin{figure}[b]
\centerline{\includegraphics[width= 8cm]{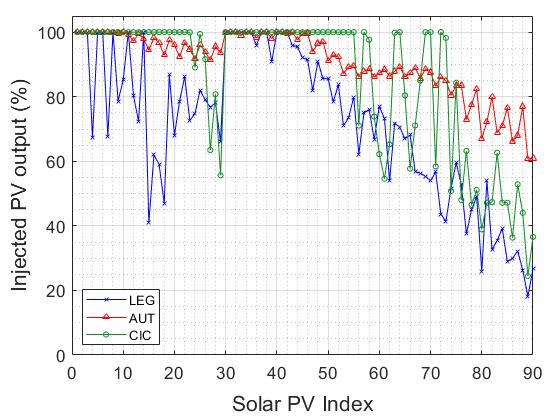}}
\caption{Results for individual customers with increasing node number roughly represents an increasing distance from the feeder level. 95mm$^2$ conductors at 80\% penetration level.} 
\label{fig:PVoutput}
\end{figure} 

All inverter control options will increase curtailment for customers further from the transformer. However, coordinated inverter control, as formulated in \eqref{eq:Objmain}, leads to a much higher variance of curtailment among customers. Some customers see up to 20\% higher curtailment than with autonomous inverters, while other customers don't experience any curtailment (Fig.~\ref{fig:PVoutput}, right). Fortunately, coordinated control has the flexibility to optimize different objectives, including objectives that encourage fairness. Such a scheme will now be investigated. 

\begin{figure}[tb]
\centerline{\includegraphics[width= 8cm]{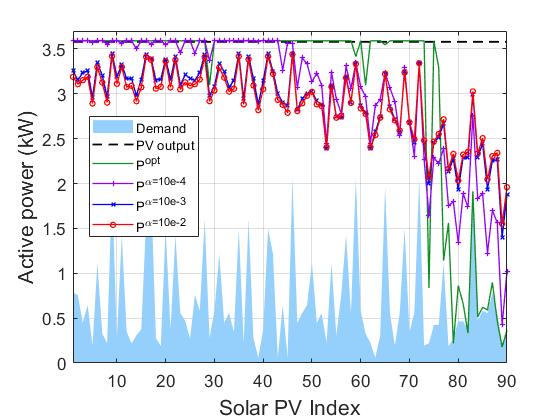}}
\caption{Active power setpoints for each inverter considering autonomous and coordinated inverter control at 80\% penetration level at 10.30\;am. $P^{opt}$ corresponds to the optimal solution for \eqref{eq:Objmain}; $P^{\alpha}$ denotes \eqref{eq:Objmain} with the fairness objective \eqref{eq:Obj4x_fair}.} 
\label{fig:nodes_fairness}
\end{figure}

The distribution of curtailment among solar PV owners for a given time instance is illustrated in Fig.~\ref{fig:nodes_fairness}. Coordinated inverter control ($P^{opt}$) targets the customers further from the distribution transformer in order to minimize line losses. A special case of coordinated inverter control is added to demonstrate that it is possible to alter the distribution of curtailment by adding a fairness objective
\begin{equation}
\nu(P_{\mathcal{C},t}) = \alpha\frac{1}{C}
\sum_{h \in \mathcal{C}} \left(\frac{P_h}{\textbf{P}^{\text{av}}_h-\textbf{P}^{\text{d}}_h} -  \frac{1}{C}\sum_{c\in\mathcal{C}} \frac{P_c}{\textbf{P}^{\text{av}}_c-\textbf{P}^{\text{d}}_c}\right)^2
\label{eq:Obj4x_fair}
\end{equation} 
to the objective of \eqref{eq:Objmain}. Its sole purpose is to redistribute curtailment across the customers with coordinated inverters reducing the lost PV output by any individual customer. The fairness objective attempts to minimize the variance of the curtailment of excess power. It is calculated as a ratio between the curtailed power $P_c$ for each PV customer $c$ and the net customer excess output ($\textbf{P}^{\text{av}}_c-\textbf{P}^{d}_c$). First, this approach ensures that all customers can meet their own demand. Second, it also rewards customers who shift their load to the periods with excess solar generation, thus reducing the potential PV curtailment.
\begin{figure}[tb]
\centerline{\includegraphics[width= 8cm]{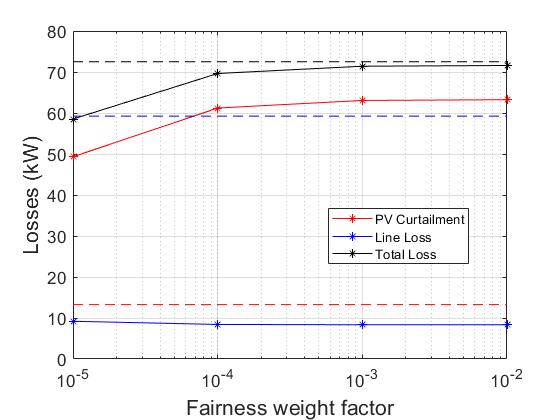}}
\caption{PV curtailment and line losses with various levels of fairness (solid lines) compared to the losses from autonomous inverter control (dashed lines).} 
\label{fig:nodes_logGraph}
\end{figure} 

A weight factor $\alpha$ controls the distribution of energy curtailed occurring along the power line with a zero value corresponding to a case without the fairness objective. A small $\alpha$ value can reduce unfair distribution by curtailing less from any individual system and redistribute the curtailment to other households, while yielding lower overall losses than autonomous inverters. The impact of different $\alpha$ values on line losses and PV curtailment is summarized in Fig.~\ref{fig:nodes_logGraph}. There is a noticeable increase in PV curtailment from the PV systems closer to the transformer to allow PV systems at the end of the line to export more power without exceeding voltage limits. When PV installation capacities vary, it is also possible to curtail to ensure either the absolute curtailment in watts is shared more equally, the generation of those curtailed is equal, or the fractional curtailment is equal. Any notion of fairness can be implemented (and traded off against performance). The goal is simply to achieve fairness towards the customers while delivering lower overall losses than currently deployed autonomous inverter control.  

The results here suggest that fair curtailment may incur high energy costs. Since the our primary concern is fair payment, simply changing the tariff structure to decouple payment from generation may be an effective solution. Because coordinated control knows the curtailment of each user, it can pay the users in proportion to the amount of electricity they \emph{would} have exported had they not been curtailed, while still generating electricity in the most efficient manner. Either way, coordination of inverter control is instrumental in ensuring fairness. 

\section{Conclusions}\label{Conclusions}
This paper investigated the benefits of coordinated inverters over autonomous inverters in terms of utilized power and voltage control. The highest reduction in overall energy losses was achieved with the smaller 70 and 95 \;mm$^2$ conductors, despite a higher curtailment with coordinated inverter control. This demonstrates the need to account for line losses and conductor size when assessing benefits of different inverter control methods.

As may be expected, coordinated inverter control is most suitable for high impedance LV networks and the LV networks constrained by the distribution transformer capacity limits. When considering a 95\;mm$^2$ overhead line, coordinated inverter control provided a modest improvement (up to 5\%) in the amount of energy fed back to the grid, relative to the autonomous Volt/VAr and Volt/Watt control. In terms of PV hosting capacity, coordinated control increased it by 20\% compared to autonomous inverters. In the case when a network without dedicated reactive power support devices is constrained by transformer capacity limits, coordinated inverters could double the PV hosting capacity.  

Coordinated inverter control also provides a guarantee that even the customers at the end of the radial network will be able to use their available output to meet their own demand. The fairness of curtailment distribution can be improved among the customers at the cost of higher overall energy losses. More importantly, coordinated control allows curtailment to be known by the distributor, meaning tariffs can be set such that optimal curtailment can be used without unfair \emph{payments} to customers. 

Future work could examine the cost of implementing coordinated inverter control and the potential savings of the CIC deployment for distribution network service providers compared to traditional network augmentation approaches for the increase of PV hosting capacity. 
\section{Appendix}\label{Appendix}
\subsection{Coordinated Inverter Control Validation}\label{Valid}
The CIC optimization algorithm was evaluated by deploying inverter target values ($P,Q$) in OpenDSS power flow software and comparing the absolute voltage difference $\Delta |V|$ and relative voltage error $\sigma$ 
in per-unit for all nodes
\begin{equation}
  \sigma = \underset{n,t}{\text{max}}\Big|\frac{V^\text{OpenDSS}_{n,t}-V_{n,t}}{V^\text{OpenDSS}_{n,t}}\Big|.
 \label{eq:sigma}
\end{equation}
The accuracy of the algorithm for balanced and unbalanced scenarios with 95\;mm$^2$ conductors is summarized in Table~\eqref{tab:valid}. We noted a tendency towards the overestimation of voltages up to a relative error of 2.1e-2 that does not have adverse effects on the coordinated inverter control at high end of the voltage range.  
%
\begin{table}[htb]
\caption{Voltage metrics for the CIC algorithm evaluation against OpenDSS power flow.}
\label{tab:valid}
\centering
\begin{tabular}{l l l l l l l}
\multicolumn{1}{c}{EU-114} &
\multicolumn{3}{c}{Balanced} &
\multicolumn{3}{c}{Unbalanced} \\
\toprule 
Pen. & 30\%& 60\% & 90\% &  30\%  & 60\% & 90\% \\
\midrule
$\sigma$ & 1.5e-2 & 1.2e-2 & 8.0e-3 & 1.9e-2 & 1.7e-2 & 1.5e-2  \\
$\text{max} \Delta |V_{k,n,t}|^+$ & 1.5e-2& 1.1e-2 & 7.7e-3& 2.1e-2 & 1.9e-2 & 1.6e-2 \\
$\text{max} \Delta |V_{k,n,t}|^-$ &6.5e-4 & 8.0e-4& 1.8e-3 & 5.3e-3 & 4.6e-3 & 4.1e-3 \\
\bottomrule
\toprule 
\multicolumn{1}{c}{IEEE-906} &
\multicolumn{3}{c}{Balanced} &
\multicolumn{3}{c}{Unbalanced} \\
\toprule 
Pen. & 30\%& 60\% & 90\% &  30\%  & 60\% & 90\% \\
\midrule
$\sigma$  & 1.3e-2& 9.9e-3& 6.4e-3 & 6.3e-3 & 6.0e-3 & 5.7e-3  \\
$\text{max} \Delta |V_{k,n,t}|^+$ &1.3e-2&9.5e-3&6.2e-3 & 3.3e-3 & 3.1e-3 & 3.1e-3  \\
$\text{max} \Delta |V_{k,n,t}|^-$ &1.7e-4 &2.9e-4 & 7.4e-4& 5.8e-3 & 5.6e-3& 5.6e-3 \\
\bottomrule
\end{tabular}
\end{table}
\subsection{Inverter reference setpoints}\label{RefSetpoints}
Inverter Volt/VAr and Volt/Watt droop curves and legacy inverter operational setpoints, given in Table \eqref{tab:active}, are in accordance to the Australian DNSPs Technical Standard  \cite{Hall2019}.
\begin{table}[htb]
\caption{Active power $P^\text{PU}$ and Reactive power $Q^\text{PU}$ reference setpoints for relevant inverter control models.}
\label{tab:active}
\centering
\begin{tabular}{l l l l l}
\toprule
  & $V$ & $\textbf{Q}^{pu}_a$  & $\textbf{P}^{pu}_a$ & $\textbf{P}^{pu}_l$\\
\midrule
$\textbf{V}^\text{min}$ &207& 0 & 0 & 0 \\
$V^\text{nom}_{t=0}$ & 230 & - & - & -\\
$\textbf{V}^\text{DB}$ &248& 0 & - & - \\
$\textbf{V}^\text{Qmin}$ &253& 0.44 & 1.0 & -\\%
$V^\text{max}_{c,t=0}$ & 257 & - & - & -\\
$\textbf{V}^\text{max}_l$ &260& - & - & 0 \\
$\textbf{V}^\text{max}_a$ &265& 0.44 & 0.2 & - \\
\midrule
$\textbf{V}^\text{trip}$ &257& 0 &  0& 0 \\
\bottomrule
\end{tabular}
\end{table}
%

\bibliographystyle{ieeetr} 
\bibliography{Coordinated.bib}
\end{document}